# EDUCATIONAL MULTI-TOUCH APPLICATIONS, NUMBER SENSE, AND THE HOMOGENIZING ROLE OF THE EDUCATOR

Anna Baccaglini-Frank

Department of Mathematics, "Sapienza" University of Rome, Italy

*As part of an educational project proposed in Italian preschools, an educator followed a tested protocol proposing two chosen iPad apps to children of ages 5-6. Though her interventions were supposedly aimed at strengthening the children's number sense, the result was a homogenization of their schemes, in various cases seemingly inhibiting development of number sense.*

## USING MULTI-TOUCH TECHNOLOGY IN THE CLASSROOM

Modern multi-touch technology offers learners new affordances that include recognition of a range of touch and multi-touch gestures as well as voice as inputs. Some studies, though not yet many, have started to analyze these affordances in relation to students' mathematical development, in particular to their development of number sense (e.g., Baccaglini-Frank & Maracci, 2015; Sinclair & Pimm, 2015; Sinclair & Baccaglini-Frank, 2016). Though it is still an elusive notion, different research communities agree that number sense is a necessary condition for learning formal arithmetic at the early elementary level and it is critical to early algebraic reasoning (English & Mulligan, 2013). In particular, literature from different fields of research converges in suggesting that using fingers for counting and representing numbers (Brissiaud, 1992), but also in more basic ways (e.g., Gracia-Bafalluy & Noël, 2008), can have a positive effect on the development of numerical abilities and of number-sense. This led us to a working hypothesis on multi-touch potential, explored in an initial study:

> "Multi-touch technology has the potential to foster important aspects of children's *development of number-sense*, including the ability to use fingers to represent numbers in an analogical format. We will call this the *multi-touch potential*." (Baccaglini-Frank& Maracci, 2015, p.6).

Our goal was to analyze the multi-touch potential of two apps for fostering preschoolers' development of number-sense, by 1) investigating the schemes that children develop in their interactions with the software, and how they use their fingers; 2) and attempting to relate the schemes enacted to the development of number-sense. The apps are *Ladybug Count* (LBC) and *Fingu* (F) (for extensive descriptions see pages 8-11 of Baccaglini-Frank & Maracci, 2015), environments providing a stimulus (either dots on the back of a ladybug or fruit in groups floating on the screen) to which the user responds placing on the screen a number of fingers that corresponds to the numerosity of the stimulus. We defined some fundamental *aspects of number sense*, and took specific ones into consideration for investigating the multi-touch potential of the two apps used. These are (see Table 1 in Baccaglini-Frank & Maracci, 2015): multiple fingers tapping (simultaneously or sequentially), subitizing (simple or double), recognizing parts of a whole, one-to-one correspondence, approximate estimation (of small or large quantities), the counting principles. In the initial study children from a class of 25, between the ages of 4 and 5 worked in groups of 5 and played 5 minutes a day for 2 weeks with the two iPad apps under the supervision of a pre-service teacher, taking turns while the other children in the group watched. Only if prompted by the pre-service teacher, a child from the group could give advice to the child playing, through verbal or gestural utterances. For this study the pre-service teacher was trained by the research group and intervened minimally during each play session.

One year later, a research-to-practice group from a university in a different city in northern Italy decided to adopt the same protocol within an educational project aimed at strengthening preschoolers' number sense, and asked for my supervision. In exchange, I was able to obtain consent forms from the parents of the children in one class to collect videos of the sessions. The 24 children of this class were in the last year of preschool (5-6 years old), from socioeconomic backgrounds comparable to those of the children in the initial study. The two major differences with respect to the initial study were: 1) the age of the children involved, 2) the fact that the educator in this study, who played the role of the pre-service teacher in the initial study, was an in-service preschool teacher (not in the same school), with a degree in psychology, and who could not avoid intervening frequently during the preschoolers' activities with the iPad. My main objective in analyzing the older children's videos was to be able to compare these children's schemes to those of the children in the initial study, and gain insight into how they evolved based on the interventions of the educator – presented to me as an "expert". I had two underlying hypotheses: 1) there would be invariant behaviors found in the enactment of (at least) the (initial) schemes of the two groups of children; 2) the educator's interventions would guide the development of certain spontaneous schemes, strengthening fundamental aspects of number sense. This paper elaborates on findings related mostly to the second hypothesis.

## THEORETICAL BACKGROUND OF THE PRESENT STUDY

I refer to the notion of *scheme* to look at and analyze the students' behavior, as in the initial study; I also add a theoretical perspective on the expected role of the teacher in mediating the children's learning processes.

### The notion of scheme

To analyze possible links between children's actions and the their goals and intentions in a given situation, and certain characteristics of the situation itself, in the initial study we used the notion of *scheme* as developed by Vergnaud (1996). Within this perspective, a scheme is seen as *an invariant organization of behavior for a given class of situations*. The main components of a scheme are: the goal and the anticipated outcomes; the rules of action, of gathering information, of control taking; and the operative invariants (implicit knowledge), including *concepts-in-actions*, that is concepts that are implicitly considered as pertinent, and *theorems-in-actions* that is, propositions believed to be true. We will refer to a (visible) recurring sequence of actions as the "enactment" of a scheme.

### Multi-touch enactments of schemes identified in the initial study

In the initial study we identified 11 enactments used by the children in LBC and 4 used in F. The ones used in LBC were classified into "general" (6), that is enactments of schemes not apparently linked to a "small" or "large" number of dots on the ladybug's back, and "specific" (5) ones that were sensitive to the number of dots to "count". This was necessary because the children seemed to hold different schemes for a very small number of dots (1-3) or large numbers of dots (7 to 10). Moreover, in several cases the children reacted to the appearance of the ladybug with a large number of dots through verbal expressions such as: "How many!" "That's a lot!". This allowed us to infer that the two situations identified above were different for them, and thus identified different schemes, possibly related to the different aspects of number-sense. For example, the most common enactment of a scheme in the presence of a small number of dots involved the rapid recognition of the small number of dots, apparently through subitizing followed by the placement of the same number of fingers (in a variety of configurations) on the screen simultaneously. This enactment did not contain verbal utterances. The most common enactment in the presence of a large number of dots involved placing all fingers on the screen and then, possibly, removing fingers one at a time until positive feedback was received from the app.

In F, selecting the exact number of fingers and placing them simultaneously on the screen within a limited amount of time was a source of difficulties for most children: indeed, it requires, among other skills, the development of advanced fine motor abilities. Typically the children's tendency was to continue to use the first enactment that seemed to be effective in a few cases, despite possible successive failures. Four enactments were identified, two of which involved recognizing the number of objects (be they in a single group or in two groups) without verbally counting, and trying to place on the screen the corresponding number of fingers either of a same hand or of two hands, simultaneously. Many children tried to place their fingers as close as possible to the floating objects, so as to "catch" them, or to reproduce with their fingers the same spatial arrangement of the floating objects (enactment 3). Only a few children quickly counted the floating objects (pointing to them and counting aloud) and then placed a same quantity of fingers on the screen (enactment 4).

**The mediating role of the teacher in students' learning processes**

In the current study my expectation was for the educator to play a fundamental, though not too invasive, role intervening to build on the children's spontaneous productions, with the aim of helping them interiorize fundamental aspects of number sense, possibly fostering thoughtful comparison between the proposed strategies. In this sense, I expected the educator to play a mediating role in the students' learning processes, picking up on the children's gestures and possibly fostering brief mathematical discussions (e.g., Bartolini Bussi, 1996, 1998).

## METHODOLOGY

All episodes in which the apps were proposed to the children were video-recorded and analyzed. The interventions of the proposing educator were flagged in the transcripts, and at the end of her work with the children, I interviewed her in order to have an additional key of interpretation.

## THE ANALYSIS OF TWO PROTOTYPICAL CASES: GIOVANNA AND SARA

The schemes enacted by Giovanna and Sara are quite similar to those of most of the other children in the class. We hypothesize this to be the case mainly because of the interventions of the educator.

**Giovanna**

Giovanna (5 years, 3 months) starts her first interaction with LBC very hesitantly, barely showing the fingers she intends to put on the screen and placing them very close together. Every time she hesitates (even for only a few seconds) the educator asks her classmates to "show Giovanna how to do this with her fingers". After three iterations of this process, Giovanna immediately looks for hints from her classmates, from one in particular, praised by the educator, checking the configuration of fingers he is raising and imitating it. Together with finger configurations, her classmates also shout out the numbers of dots. Giovanna pretends to count them, pointing to a few and then repeats the number pronounced by her classmates. The educator says nothing to stop the classmates from talking, glances at the screen to check it, and makes comments like:

| Educator: | Come on, you know how to do nine [with your fingers]! |

The educator also highly praises children who "know", as in the following example.

| | [A ladybug with 7 dots on her back appears] |
| Giovanna: | [Counts up the dots pointing with her finger.] Seven [in a whisper]. |
| Teacher: | Seven! How do you do seven with your fingers, Giovanna? Come on! |
| Giovanna: | [She timidly raises the fingers on a hand and the index and middle finger of her other hand, imitating the fingers shown by a classmate.] |
| Teacher: | Right, very good! See, you *know*! |

These interventions by the educator both in LBC and in F seem to lead Giovanna to develop two schemes described in the box below, with the goal of doing what she thinks the educator (seen as a teacher) wants, and seemingly identifying two situations in LBC, based on whether she immediately recognizes the number of dots on the ladybug's back (1a) or not (1b).

1. Figure out the number of dots/fruits, to do this either
      1a. recognize the number immediately
      1b. count them up from one;
2. say the number,
3. use *the* fixed configuration approved by the educator for that number.

She seems to also develop what we could see as a concept-in-action:

Every number pronounced verbally corresponds to a *fixed* configuration of fingers.

Giovanna seems to inhibit any enactment that involves configurations of fingers other than what she believes to be *the right one*. In fact, in one episode during her first playing session, Giovanna sees the dots on the ladybug's back [there are two on each wing], and she timidly raises two fingers on each hand fingers, waiting for feedback. She sees other children showing four fingers raised on one hand, and copies their finger configuration. The educator says nothing to the other children and simply says "Good!" to Giovanna when she touches the screen and receives positive feedback. A short time after, a ladybug with 9 dots appears. Giovanna does not count, but she over estimates, placing all fingers of both hands on the screen. Instead of trying to adjust her fingers (as in schemes identified in the initial study), for example, by lifting one, she takes her hands off and looks for a configuration to copy. After these two episodes Giovanna's behavior can be described entirely in light of the schemes we hypothesized above.

In general, when Giovanna cannot remember what she believes to be *the approve*d finger configuration, she depends on her classmates' hints, and copies, seemingly with no control over the answer that she then gives; or, if she cannot catch a hint quickly, she listens to the number pronounced by the class and counts up from "one", raising her fingers one at a time and always in the same order. Clearly, remembering fixed configurations of fingers associated to a word requires a lot of (otherwise unnecessary) memory and it could even inhibit the development of fundamental aspects of number sense: the child can be successful without ever putting in relation the dots and the fingers raised, other than through a verbal utterance (when a finger configuration is not shown directly), and without developing, for example, awareness of part-whole relations.

In F, Giovanna relies entirely on reaching (either by herself or hearing it from her classmates) a verbal pronunciation of the number of floating fruits from which she produces *the* configuration of fingers if she remembers it. She does not count up her fingers because the game does not give her the time. Her configurations of fingers have nothing to do with the partitions of the fruits into smaller sets when there is more than one. For example, when two floating sets of 2 fruits appear, she hears a classmate say: "four" and puts down the fingers of her right hand excluding her thumb. The same happens when sets of 3 and 1 appear on the screen. Giovanna seems to rely heavily on her schemes; so much that when possibly perturbing episodes occur, the scheme remains unchanged. For example, when 4 and 1 fruits appear she hears a classmate say (erroneously) "four" and she puts down her usual "four" configuration. When 3 and 2 fruits appear a classmate shows 3 fingers on one hand and 2 on the other; Giovanna sees, but she hesitates and then says: "There are five" and places down all her fingers of a single hand.

Her schemes turn out to be successful in F and Giovanna passes to level II of the game. Now more than 5 fruits can appear. The other children no longer have time to figure out how many fruits are floating around on the screen before they disappear, so Giovanna counts them up each time, starting from "one" and saying the numbers aloud as she points to each fruit. She appears to not know (at least not quickly enough) the configurations for any of the quantities above 5, so she either puts down no fingers or she tries to put down some, frequently less than 5 and loses.

In F her enactments do not seem to change, however in LBC they do. During the following playing sessions, Giovanna proceeds more and more independently with respect to what her classmates say, and the form of her gestures change, as well. For example, she places her fingers on the screen spreading them out on the leaf rather than bunching them up, like during the initial episodes. During the last session with LBC Giovanna has learned to generate fixed configurations quickly for quantities below 5, and for larger quantities she modifies her behavior for figuring out how many fingers to raise. Interestingly, she never counts up her fingers when she recognizes the number of dots. These changes can be seen in terms of new schemes (1a-2a; 1b-2a; 1b-2b):

1. Figure out the number of dots, to do this either

      1a. recognize the number immediately (if below five) and say it out loud;

      1b. (if 1a fails) count them up from one;

2a. use *the* fixed configuration of fingers, when known, corresponding to the number pronounced,

2b. otherwise count up fingers (in a constant order) starting from one.

For example, when a ladybug with 8 dots appears, Giovanna counts up the dots and counts her fingers, starting with "one" as the thumb of her right hand. The educator gives very positive feedback. The same happens when a ladybug with 10 dots appears, and then when a ladybug with 8 dots appears again. When the next ladybug appears and it presents, again, 8 dots, Giovanna seems to recognize the configuration and remember the configuration for "eight" (as a hand and 3 fingers). This suggests that she still holds valid the concept-in-action we had hypothesized earlier.

This can also be seen in how Giovanna seems to privilege schemes 1a-2a and 1b-2a over ones using 2b. For example, when a ladybug with 4 dots appears she says "four" and in a seemingly automatic way counts up four fingers starting from her thumb on the right hand, but as soon as she sees the four fingers she changes the fingers to her *fixed configuration* for "four".

**Sara**

Although Sara (5 years, 5 months) is not in Giovanna's group, with her the educator keeps on intervening in the same way as with Giovanna, proposing to count the dots or fingers immediately at the smallest hesitation, and praising her emphatically whenever she does count. During her initial interactions both with LBC and F, Sara is less insecure than Giovanna: for quantities of 1, 2, 3 or 4 she simply says aloud the number corresponding to the quantity and raises a known configuration of fingers. Sara uses constant configurations for "one", "two" and "three", while for "four" she seems to flexibly change the fingers raised and placed on the screen. The educator, in these cases, simply praises Sara for getting positive feedback from the software.

Both for LBC and F the schemes developed by Sara seem to be very similar to Giovanna's:

| |
|---|
| 1. Figure out the number of dots/fruits, to do this either |
|        1a. recognize the number immediately (if below five) and say it out loud; |
|        1b. (if 1a fails) count them up from one, aloud, pointing to each; |
| 2a. raise *any* known configuration of fingers corresponding to the number pronounced, |
| 2b. otherwise count up fingers (in a constant order) starting from one, |
| 3. In any case, count up fingers before placing any on the screen (if in 2a, this is to please educator). |

Although Sara shows a bit more flexibility than Giovanna in making appropriate finger configurations for quantities up to 4, she, too, seems to be conditioned by the educator's insistence on counting. This results in episodes such as the following. In LBC a ladybug with 4 dots appears:

| | |
|---|---|
| Sara: | Four [and raises two fingers (index and middle) on both hands immediately]. One, two, three, four [she counts the fingers on one hand starting from the thumb]. |
| Educator: | Very good! |
| Sara: | [She switches back to her initial configuration of 2 and 2 fingers and places them on the screen, receiving positive feedback from LBC]. |
| Educator: | Oh, that's OK, too. Good job! |

All the counting in this enactment is too time consuming to be effective in F, so Sara receives negative feedback almost every time. She soon asks to stop playing and to give another classmate a turn. The educator satisfies her request and calls a classmate to play.

**The other cases**

Various children's initial schemes have enactments similar to those found in the first study, but these are quickly adapted to schemes like those of Giovanna or Sara, so that they include counting dots/fruits and, usually, fingers before placing them on the screen simultaneously. However, four children simply refuse to count, possibly because they have not sufficiently mastered the counting principles. The children either imitate the configurations of fingers shown by their classmates, when they were able to, and failed otherwise, or they continued to use enactments similar to those of children in the initial study. For example, Amanda (5 years, 2 months), never says numbers aloud, but places down precise numbers of fingers, in a variety of configurations, for quantities up to 4, both in LBC and in F, while for larger numbers she seems to estimate, quickly putting down a hand of fingers and some additional ones. Whenever this fails, the educator tells her to count (this happens especially in LBC), but she continues not to, and instead tries to recognize recurring configurations of dots and remember and use the configuration of fingers that worked previously.

## THE INTERVENTIONS OF THE EDUCATOR

In general, the educator's behavior can be described as follows. She never takes the time to discuss the children's schemes or enactments, neither collectively, nor individually. Moreover, although she (fortunately!) accepts different finger configurations for a same quantity, she does not ever explicitly comment on how a same quantity can be represented through different finger configurations, for example, putting in one-to-one correspondence two different quantities of fingers for a given quantity of dots or fruits. The only "sharing" that the educator fosters is in cases in which the child playing appears to be hesitant: she calls on other children to "show your classmate how to do it". From the follow-up interview we found that the educator did have a particular procedure in mind that she thought was best, because applicable to all the situations generated by the apps. Below are various claims she made in the interview.

> Educator: You need to count up the dots or fruits and then count up your fingers, quickly, and put them down together. […] Once a kid knows how to do a certain number with his fingers, he doesn't have to count up his fingers any more, and he can use whatever way he wants. […] which fingers they raise is not important, but they should do it quickly and there are easier ways of doing it. […] This is how they can always experience success.

## DISCUSSION AND CONCLUSION

Although the educator was a trained psychologist and had been introduced as an "experienced" preschool teacher, the hypothesis that her interventions would promote certain enacted schemes, strengthening fundamental aspects of number sense, possibly through brief mathematical discussions (Bartolini Bussi, 1998) turned out to be rather naïve and excessively optimistic. In fact, the children did initially enact schemes quite similar to those of the children in the first study, however the educator did not pick up on most of these enactments, and vigorously promoted *counting* and/or use of *known finger configurations*, overruling the children's strategies. Her interventions were mostly consistent with the claims she made during the follow-up interview. Her main goal was to help children experience success in the apps and, with respect to strengthening number sense she intended to help children learn to count and to represent numbers on their fingers. She did not mention, for example, recognition of part-whole relationships, double subitizing, estimating, or even counting on; and indeed she fostered none of these.

Her short-term goal of helping the children experience success, and her narrow-sighted view of how to obtain this while fostering the development of number sense, actually slowed down or inhibited altogether such development (at least during this experience) for many children. In most cases the children seemed to develop schemes like those of Sara and Giovanna, whose enactments included trying to memorize fixed configurations of fingers corresponding to verbal numbers (this was a common children's misinterpretation of various of the educator's comments), or simply copying. The other abilities the educator explicitly intended to promote were to represent a same number in different ways (with fingers), and to count. She failed to promote the former by not taking the time to discuss any of the situations in which children used different representations of a same number. As for counting, one might say that at least this ability was promoted and most children learned to count faster. Though this is true, it is not clear how such ability enhances number sense: indeed it is important to be able to count, but achieving mastery in counting does not simply mean learning to do it *faster*! For example, it is also important to learn to count on from a number greater than "one", to count backwards, and to learn to r*eplace* counting with more effective strategies (e.g., Gray & Tall, 1994).

Though the educator claimed that she valued sharing and discussing children's strategies, not in a single episode could she be seen doing this (at least not in the sense of Bartolini Bussi, 1998), and so her mediation ended up being quite different than what I expected. A number of questions about why the educator *did* and *did not* do certain things come to mind, but above all, as math educators we probably agree that homogenizing children's solution procedures as was done here is not what education should seek for. However, the phenomenon does seem representative of what frequently happens also at higher levels of mathematics teaching and learning. I found it quite sad to discover how strong the phenomenon can be even before formal schooling starts, even when an educational intervention is carried out by someone socially described as "expert". Finally, I dare not ask what might happen to this educator (and the children she works with) in a more open digital environment where many different tasks can be proposed and a greater variety of solutions can be given.